\documentclass[12pt]{article}

\usepackage[T1,T2A]{fontenc}
\usepackage[utf8]{inputenc}

\usepackage{amsfonts}
\usepackage{amssymb}
\usepackage{amsmath}
\usepackage{amsthm}

\def \le {\leqslant}
\def \ge {\geqslant}

\usepackage{tikz} 
\usetikzlibrary{calc,arrows,decorations.pathreplacing,fadings,3d,positioning}

\begin{document}

 \begin{Large}

 \centerline{\bf \"Uber die  Funktionen des Irrationalit\"atsma\ss es}
  \centerline{\bf f\"{u}r zwei irrationalen Zahlen}
 \vskip+0.5cm

 \centerline{von Nikolay Moshchevitin\footnote
 {Steklow-Institut f\"ur Mathematik der Russischen Akademie der Wissenschaften. Diese Arbeit wurde unterst\"utzt durch RNF Grant No. 14-11-00433.
 }}
 
 \vskip+1cm
 \begin{small}
 {\bf Abstract:}\,
 For real $\xi$ we consider irrationality measure function
 $\psi_\xi (t) = \min_{1\le q \le t, \, q\in \mathbb{Z}} ||q\xi||$.
 We prove that  in the case $\alpha \pm \beta \not\in \mathbb{Z}$ there exist arbitrary large  values of $t$ with
 $|\psi_\alpha (t) -\psi_\beta(t)| \ge \left(\sqrt{\frac{\sqrt{5}+1}{2}}-1\right) \min (\psi_\alpha(t), \psi_\beta (t))$. This result is optimal.
  
 \end{small}
\vskip+1cm 
 \end{Large}

    {\bf 1. Resultate.}
    
    F\"{u}r eine irrationale Zahl $\xi \in \mathbb{R}$  betrachten wir die  Funktion des Irrationalit\"atsma\ss es
 $$
 \psi_\xi (t) = \min_{1\le q \le t, \, q\in \mathbb{Z}} ||q\xi||
 .
 $$
 Sei
 $$
 Q_0\le Q_1<Q_2<...< Q_n < Q_{n+1}<...
 $$
 die Folge der Nenner der Naherungsbr\"{u}che  f\"{u}r $\xi$. 
  Aus dem Lagrangeschen Gesetz der besten N\"aherungen (siehe zum Beispiel \cite{P} oder \cite{RS})
folgt, dass
$$
\psi_\xi (t) = ||Q_n \xi ||\,\,\,\text{f\"ur}\,\,\, Q_n\le t < Q_{n+1}.
$$
 Aus den
 Minkowskischen Gitterpunktsatz 
 kennen wir, dass
  f\"{u}r jeden $\xi\in \mathbb{R}$ und 
    f\"{u}r jeden $t\ge 1$ gilt 
 \begin{equation}\label{0}
 \psi_\xi (t) \le \frac{1}{t}.
 \end{equation}
 
  2010 hatten Kan und Moshchevitin \cite{KM}  die folgende Behauptung bewiesen.
  
  {\bf Satz A.}
 \,{
 \it Seien $ \alpha$ und $\beta$ zwei irrationale Zahlen. Dann falls $ \alpha\pm \beta \not\in \mathbb{Z}$,
die Differenz
 $$
 \psi_\alpha (t) -\psi_\beta (t)
 $$
 wechselt 
 ihr Vorzeichen  
 unendlich viele Male wenn $t \to \infty$.}
 
 Verschiedene metrische Ergebnisse  \"{u}ber die Funktion  $\psi_\alpha (t) $ wurden in \cite{Den} betrachten.
 
 2017 gab Dubickas \cite{D} einen anderen Beweis des Satzes A  mit der Kombinatorik der Worten und hat das folgenden Ergebnis bewiesen.
 
  {\bf Satz B.}
 \,
{ \it Seien $ \alpha$ und $\beta$ zwei irrationale Zahlen. Dann falls $ \alpha\pm \beta \not\in \mathbb{Z}$,
 hat man
 $$
 \limsup_{t\to+\infty}\left|
 \frac{1}{\psi_\alpha (t)} -\frac{1}{\psi_\beta (t)}\right| = +\infty
 $$
 }

  Hier beweisen wir eine starkere Behauptung.
 Betrachten wir die  Zahlen
  \begin{equation}\label{KA}
  \tau =  \frac{\sqrt{5}+1}{2}\,\,\,\,\,
  \text{und}\,\,\,\,\,
  K =  \sqrt{\tau} - 1 =\sqrt{\frac{\sqrt{5}+1}{2}}-1= 0.2720^+.
  \end{equation}
  
  {\bf Satz 1.} \,{\it 
  Seien $ \alpha$ und $\beta$ zwei irrationale Zahlen. Dann falls $ \alpha\pm \beta \not\in \mathbb{Z}$,
  zu jedem $T \ge 1$ gibt es  $t\ge T$ mit
 \begin{equation}\label{Au}
  |
  \psi_\alpha (t) -\psi_\beta (t)|
  \ge K\cdot \min (\psi_\alpha(t), \psi_\beta (t)).
  \end{equation}
  }
  
  Aus Satz  1 und (\ref{0}) ableiten wir 
  
  {\bf Folgerung.}\,
  {\it
  Seien $ \alpha$ und $\beta$ zwei irrationale Zahlen. Dann falls $ \alpha\pm \beta \not\in \mathbb{Z}$,
  zu jedem $T \ge 1$ gibt es  $t\ge T$ mit
   \begin{equation}\label{B}
  \left|\frac{1}{
  \psi_\alpha (t)} -
  \frac{1}{\psi_\beta (t)}\right|
  \ge {K}{t}.
  \end{equation}
  }

 Das
Ergebnis des Satzes 1  ist optimal.

  {\bf Satz 2.} \,{\it In (\ref{Au}) kann $K$ durch keine gro\ss ere Zahl ersetzt werden.  }

    {\bf 2. Kettenbr\"{u}che und Hilfss\"{a}tze.}
    Betrachten wir die Kettenbr\"{u}che
    $$
    \alpha=[a_0;a_1,a_2,a_3,...],\,\,\,\,\,\,\,\,  \,\,\,\,
    \beta = [b_0;b_1,b_2,b_3,...]
    $$
    und die Naherungsbr\"{u}che 
    $$
     \frac{p_n}{q_n}=[a_0;a_1,a_2,a_3,..., a_n],\,\,\,\,\,\,\,\, \,\,\,\,
    \frac{t_m}{s_m} = [b_0;b_1,b_2,b_3,..., b_m]
    .
    $$
    Definieren wir
    $$
      \alpha_n=[a_n;a_{n+1},a_{n+2},a_{n+3},...],\,\,\,\,\,\,\,\,  \,\,\,\,
    \beta_m = [b_m;b_{m+1},b_{m+2},b_{m+3},...],
    $$
    $$
    \xi_n = |q_n\alpha - p_n|
      ,\,\,\,\,\,\,\,\,  \,\,\,\,
    \eta_m = |s_m \beta - r_m|.
    $$
    Es ist klar, dass
      \begin{equation}\label{kkl}
      \frac{\xi_{n-1}}{\xi_n} =\alpha_{n+1},
      \,\,\,\,\,\,\,\,  \,\,\,\,
        \frac{\xi_{m-1}}{\xi_m} =\beta_{m+1}.
      \end{equation}
    
Betrachten wir die Werte   
  \begin{equation}\label{1}
 \lambda_n (\alpha) = 
 \sqrt{\alpha_{n+1}}-1,
 \,\,\,\,\,\,\,\,  \,\,\,\,
  \lambda_m(\beta) = 
 \sqrt{\beta_{m+1}}-1.
  \end{equation}

   \vskip+0.3cm
  {\bf Hilfssatz 1.}
  {\it  
  
  \noindent
  {\rm({\bf i})}
  F\"{u}r
  jeden  $\eta \ge \xi_{\nu-1}$ hat man
  $$
 \eta- \xi_{\nu} > \lambda_\nu (\alpha) \xi_\nu;
 $$
 
   \noindent
  {\rm ({\bf ii})}
  f\"{u}r
  jeden $\eta \le\xi_{\nu}$ hat man
  $$
 \xi_{\nu-1} -\eta> \lambda_\nu (\alpha) \eta;
 $$
   \noindent
  {\rm ({\bf iii})}
  f\"{u}r
  jeden  $\eta \ \in (\xi_\nu,\xi_{\nu-1})$ hat man
  $$
  \min\left(
  \frac{
 \xi_{\nu-1} -\eta}{\eta},
  \frac{
 \eta-\xi_\nu}{\xi_\nu}\right)
\ge \lambda_\nu (\alpha) .
 $$ 
    }

   \vskip+0.3cm

  Beweis. \, 
  Die Behauptungen
   {\rm ({\bf i})} und
   {\rm ({\bf ii})}
   sind klar als
   $ \xi_{\nu-1} - \xi_\nu =(\alpha_{\nu+1} -1)\xi_\nu \ge \lambda_\nu (\alpha) \xi_\nu
   $.

  Bemerken wir, dass
  $$
  \frac{\xi_{\nu-1}}{\eta}\cdot \frac{\eta}{\xi_\nu} = 
  \frac{\xi_{\nu-1}}{\xi_\nu} = \alpha_{\nu+1}.
  $$
  Dann
  f\"{u}r    $\eta \ \in (\xi_\nu,\xi_{\nu-1})$ 
  haben wir entweder
  $\frac{\xi_{\nu-1}}{\eta} \ge \sqrt{\alpha_{\nu+1}}$ und
  $\xi_{\nu-1}-\eta \ge \lambda_n (\alpha)  \eta$, oder
   $\frac{\eta}{\xi_\nu} \ge \sqrt{\alpha_{\nu+1}}$ und
  $\eta-\xi_\nu \ge \lambda_n (\alpha)  \xi_\nu$.$\Box$

     \vskip+0.3cm
   
     {\bf Folgerung 1.}\,\,
  {\it Sei $\lambda_n (\alpha) \ge K$.
  Falls $q_n$  ist kein Nenner des Naherungsbruchs f\"{u}r $\beta$,
   hat man 
   \begin{equation}\label{enn}
   \begin{array}{c}
   \text{
   entweder
  }\,\,\,\,\,
     |\psi_\alpha(t) -\psi_\beta (t)| \ge K\cdot  \min (\psi_\alpha(t), \psi_\beta (t))\,\,\,\,
  \text{f\"{u}r}\,\,\,\,  t \in [q_n-1, q_n),\cr
  \text{
  oder
  }\,\,\,\,
|\psi_\alpha(q_n) -\psi_\beta(q_n)| \ge K\cdot   \min (\psi_\alpha(q_n) , \psi_\beta (q_n))
 .\,\,\,\,\,\,\,\,\,\,\,\,\,\,\,\,\,\,\,\,\,\,\,\,\,\,\,\,\,\,\,\,\,\,\,\,\,\,\,\,\,\,\,\,\,\,\,\,\,\,\,\,
 \end{array}
  \end{equation}
  }
  
   \vskip+0.3cm

  Nat\"{u}rlich, gibt es eine 
  \"{a}hnlich Behauptung  \"{u}ber $\beta$ statt $\alpha$.

   \vskip+0.3cm
  
       {\bf Folgerung 2.}\,\,
  {\it Sei $\lambda_m (\beta) \ge K$.
  Falls $s_m$  ist kein Nenner des Naherungsbruchs f\"{u}r $\alpha$,
   hat man 
   \begin{equation}\label{ennq}
   \begin{array}{c}
   \text{
   entweder
  }\,\,\,\,\,
     |\psi_\beta(t) -\psi_\alpha(t)| \ge K\cdot  \min (\psi_\alpha(t), \psi_\beta (t))\,\,\,\,
  \text{f\"{u}r}\,\,\,\,  t \in [s_m-1, s_m),\cr
  \text{
  oder
  }\,\,\,\,
|\psi_\beta(s_m) -\psi_\alpha(s_m)| \ge K\cdot  \min (\psi_\alpha(s_m), \psi_\beta (s_m))
 .\,\,\,\,\,\,\,\,\,\,\,\,\,\,\,\,\,\,\,\,\,\,\,\,\,\,\,\,\,\,\,\,\,\,\,\,\,\,\,\,\,\,\,\,\,\,\,\,\,
 \end{array}
  \end{equation}
  }
   \vskip+0.3cm
  
   Beweis der Folgerung 1.
   Weil $q_n$   kein Nenner des Naherungsbruchs f\"{u}r $\beta$ ist,
   ist die Funktion $\psi_{\beta} (t)$ f\"{u}r $ t \in [q_n-1, q_n]$ konstant, sodass
   $$
    \psi_{\beta} (t) = \psi_\beta(q_n),\,\,\,\,\,t \in [q_n-1, q_n].
    $$
    Es ist klar, dass
   $$
   \psi_\alpha(t) = \xi_{\nu-1},\,\,\,\,\,  t \in [q_n-1, q_n);\,\,\,\,\,\,\,\,\,\,\, \psi_\alpha (q_n) = \xi_n.
   $$ 
   Nehmen wir $\eta = \psi_\beta(q_n)$. Die Aussage folgt daraus.$\Box$

   \vskip+0.3cm
  {\bf Hilfssatz 2.}

  \noindent
  ({\bf i})
  {\it  Sei}
  $ \lambda_n  (\alpha) < \sqrt{2}-1$. {\it Dann ist} $ a_{n+1} = 1$.
    
    \noindent
    ({\bf ii}) {\it F\"{u}r jeden $n$ gilt 
  $\max(\lambda_{n-1} (\alpha) ,\lambda_{n} (\alpha) ) \ge K$.
  }
  \vskip+0.3cm
  
  Beweis.\,
    ({\bf i}) folgt  aus der  letzten Gleichung aus (\ref{1}). 
    Falls $\lambda_{n-1} (\alpha)  \ge \sqrt{2}-1 $ ist ({\bf ii}) klar als $ \sqrt{2}-1 > K$.
    Um    ({\bf ii}) zu beweisen f\"{u}r $\lambda_{n-1} < \sqrt{2}-1 $, bemerken wir, dass $ a_{n} = 1$ nach   ({\bf i}). Dann hat man
    $\alpha_n = 1+\frac{1}{\alpha_{n+1}}$ und
    $
    \max(\alpha_{n},\alpha_{n+1}) \ge \tau.
    $
  Alles ist bewiesen. $\Box$

     \vskip+0.3cm
  {\bf Hilfssatz 3.}
  
  {\it Betrachten wir die Kettenbr\"{u}che
  $$
  \frac{p_n}{q_n} = [0; a_1,...a_{n-2},a_{n-1}, a_n],\,\,\,
    \frac{p_{n-1}}{q_{n-1}} = [0; a_1,...,a_{n-2},a_{n-1}],\,\,\,
    \frac{p_{n-2}}{q_{n-2}} = [0; a_1,...,a_{n-2}]
  $$
  und
    $$
  \frac{r_m}{s_m} = [0; b_1,...,b_{m-1}, b_m],\,\,\,\,\,
    \frac{r_{m-1}}{s_{m-1}} = [0; b_1,...,b_{m-1}].
  $$}
  
  \noindent
      ({\bf i}) {\it
  Falls $ q_n = s_m$ und $ q_{n-1} = s_{m-1}$, hat man
  
  entweder $ m=n$ und $ a_j= b_j,\,\,\, 1\le j \le n$,
  
  oder $ m=n+1, a_1\ge 2$  und $
 ( b_1, b_2,  ...  ,b_{n+1}) = (1, a_1-1, ..., a_n)$,
 
 oder
  $ m+1=n, b_1\ge 2$  und $
 ( 1,b_1-1,  ...  ,b_{m}) = (a_1,a_2 ..., a_{m+1})$.
 }
 
   \noindent
  ({\bf ii}) {\it
  Falls $a_n =1$ und $ q_n =s_m, q_{n-2} = s_{m-1}$,  hat man
  
    entweder $ m=n$, $ a_1 \ge 2$ und  $( b_1, b_2,b_3,  ...  ,b_{n-1},b_n) = (1, a_1-1,a_2, ...,a_{n-2} , a_{n-1}+1)$,
  
  oder $ m=n-1, a_1\ge 2$  und $
 ( b_1,  ...  ,b_{n -2},b_{n-1}) = ( a_1, ...,a_{n-2}, a_{n-1}+1)$,
 
 oder
  $ m=n-2, b_1\ge 2$  und $
 ( 1,b_1-1,  ...  ,b_{m-1},b_{m}) = (a_1,a_2 ..., a_{m} a_{m+1}+1)$.
 }
      \vskip+0.3cm
      
      Beweis.\,  ({\bf i}) folgt aus den Gleichungen
      $$
      \frac{q_{n-1}}{q_n} = [0;a_n,a_{n-1},...,a_1] =
        \frac{s_{m-1}}{s_m} = [0;b_m,b_{m-1},...,b_1].  
      $$
      Um ({\bf ii}) zu beweisen, soll man bemerken, dass
      $$
      \frac{q_{n-2}}{q_n} = \frac{q_{n-2}}{q_{n-1}+q_{n-2}}=
      \frac{1}{1+q_{n-1}/q_{n-2}}
       = [0; a_{n-1}+1,a_{n-2},...,a_1] =\frac{s_{m-1}}{s_m} = [0;b_m,b_{m-1},...,b_1].\,\,
       \Box
      $$
  
   {\bf 3. Beweis des Satzes 1.}

   Wir
   analysieren den Kettenbruch f\"{ur} $\alpha$ 
   mit die Werte $\lambda_\nu = \lambda_\nu (\alpha)$ zusammen
   und
    betrachten verschiedene Falle.

   Es folgt nach dem Hilfssatz 2, dass 
   entweder
   
   \noindent
   ({\bf Fall 1})
   gibt es  unendlich viele $\nu$ mit  $\min(\lambda_{\nu-1},\lambda_\nu) \ge K$,

   oder

      \noindent
      ({\bf Fall 2})
   gibt es $ j_0$ mit den folgenden Eigenschaften
   
      \noindent
   1) $\lambda_{j_0+2j }\ge K$ f\"{u}r alle $ j =0,1,2,3,...$,
   
      \noindent
   2) $ a_{j_0+2j+2} = 1 $ f\"{u}r alle $ j =0,1,2,3,...$.

   Betrachten wir {\bf Fall 1}.

     \noindent
   {\bf Teilfall 1.1.} 
Es gibt unendlich viele $\nu$ mit
$\min(\lambda_{\nu-1},\lambda_\nu) \ge K$
und   es gibt $l\in \{ \nu-1, \nu\}$ 
dass
$q_n $ ist kein Nenner des Naherungsbruchs f\"{u}r $\beta$. 
 Dann nach der Folgerung 1 
 hat man (\ref{enn}) mit $ n=l$.

   \noindent
   {\bf Teilfall 1.2.} 
 F\"{u}r jeden diesen $\nu$  sind die  beide Nenner $q_{\nu-1}$ und $q_\nu$ auch die Nenner des Naherungsbruchs f\"{u}r $\beta$.
 Dann  gibt es $\mu = \mu(\nu)<\kappa =\kappa (\nu)$ mit
 $$
 q_{\nu-1} = s_\mu,\,\,\,\,\,
 q_{\nu}=s_\kappa.
 $$
 
 \noindent  
{\bf Teilfall 1.2.1.} 
Es gibt unendlich viele $\nu$ mit $\kappa (\nu) = \mu(\nu) + 1$.
In diesem Fall nach  dem Hilfssatz 3 ({\bf i}) hat man
$\alpha \pm \beta \in \mathbb{Z}$.

\noindent  
{\bf Teilfall 1.2.2.} 
Es gibt unendlich viele $\nu$ mit $\kappa (\nu) =\mu(\nu) + 2$.
Dann gibt es $\sigma\in \mathbb{Z}, \mu < \sigma< \kappa$.
Falls $b_{\sigma+1} \ge 2$, alles ist klar.
Wir benutzen 
den Hilfssatz 1 ({\bf i}) und der Folgerung 2, um (\ref{ennq})
zu erhalten. 
Falls $b_{\sigma+1} =1$,
 wir benutzen den Hilfssatz 3 mit $\alpha$ statt $\beta$ und $\beta$ statt $\alpha$,
 um $\alpha \pm \beta \in \mathbb{Z}$ zu sehen.

  \noindent  
{\bf Teilfall 1.2.3.} 
Es gibt unendlich viele $\nu$ mit $\kappa (\nu) \ge\mu(\nu) + 3$.
Dann gibt es $\sigma\in \mathbb{Z}, \mu < \sigma< \kappa-1$.
Nach dem Hilfssatz 2 ({\bf ii}) f\"{u}r $\beta$ statt $\alpha$ und der Folgerung  2 
mit $ m = \sigma$
hat man  (\ref{ennq}).

Nun 
   betrachten wir {\bf Fall 2}.
   
     \noindent
   {\bf Teilfall 2.1.} 
Es gibt unendlich viele $j$, sodass $q_{j_0+2j}$
kein Nenner des Naherungsbruchs f\"{u}r $\beta$ ist.
Wegen $\lambda_{j_0+2j}\ge K$,
kann man  die  Folgerung 1 anwenden.
So  
 gibt es (\ref{enn}) mit $ n = j_0+2j$
 f\"{u}r diese $j$.
 
  \noindent
   {\bf Teilfall 2.2.} 
   F\"{u}r alle $j$, sind $q_{j_0+2j}$  die Nenner der Naherungsbr\"{u}che f\"{u}r $\beta$.
   Dann m\"{u}ssen wir mehrere Falle betrachten.

  \noindent
   {\bf Teilfall 2.2.1.}
   Es gibt unendlich viele $j$ mit
   $$
   q_{j_0 +2j} = s_\mu, \,\,\,\,\,
      q_{j_0 +2j+2} = s_{\mu+1}.
      $$
      Im Fall 2 haben wir $ a_{j_0+2j+2} = 1$. Nun aus  dem Hilfssatz 3 ({\bf ii}) hat man
      $\alpha\pm \beta \in \mathbb{Z}$.
      
        \noindent
   {\bf Teilfall 2.2.2.}
   Es gibt unendlich viele $j$ mit
   $$
   q_{j_0 +2j} = s_\mu, \,\,\,\,\,
      q_{j_0 +2j+2} = s_{\mu+2}.
      $$
      Wegen $ a_{j_0+2j+2} =1$, haben wir
      $$
       s_{\mu+2} = b_{\mu+2}s_{\mu+1}+q_{j_0 +2j} = 
      q_{j_0 +2j+2} =   q_{j_0 +2j+1}+  q_{j_0 +2j}.
      $$
      Falls $ b_{\mu+2}= 1$ f\"{u}r unendlich viele $j$,
      haben wir  $\alpha\pm \beta \in \mathbb{Z}$
      nach dem  Hilfssatz 3 ({\bf i}).
      Falls $ b_{\mu+2}\ge 2$ hat man
      $q_{j_0+2j}< s_{\mu+1}< q_{j_0+2j+1}$ und nach die Folgerung 2 
      gilt (\ref{ennq}) mit $ m = \mu+1$.

            \noindent
   {\bf Teilfall 2.2.3.}
   Es gibt unendlich viele $j$ mit
   $$
   q_{j_0 +2j} = s_\mu, \,\,\,\,\,
      q_{j_0 +2j+2} = s_{\mu+3},\,\,\, \mu = \mu (j).
      $$
      Falls $
   q_{j_0 +2j+1} = s_{\mu+1}$ f\"{u}r unendlich viele $j$,
   haben wir entweder $ b_{\mu+1} = 1$ und
   nach dem  Hilfssatz 3 ({\bf i}) f\"{u}r $\beta$ statt $\alpha$ hat man $ \alpha \pm \beta = 1$,
   oder $ b_{\mu+1} \ge 2$, und dann (\ref{ennq}) folgt.
   Falls  
   $$
   q_{j_0 +2j} = s_\mu <s_{\mu+1}< s_{\mu+2}<
   q_{j_0+2j+1},\,\,\,\,\,
   \text{oder}\,\,\,\,\,
         q_{j_0 +2j+1} < s_\mu <s_{\mu+1}< s_{\mu+2}<q_{j_0 +2j+2} = s_{\mu+3},
         $$
         aus den Hilfss\"{a}tzen 1,2 finden wir einen Index $m \in\{\mu+1,\mu+2\}$
         f\"{u}r den  (\ref{ennq}) gilt.
         Betrachten wir den letzten Teilfall wo f\"{u}r  unendlich viele $j$  
           $$
   q_{j_0 +2j} = s_\mu <s_{\mu+1} <
   q_{j_0+2j+1}< s_{\mu+2} <q_{j_0 +2j+2} = s_{\mu+3}
         $$
        mit $ \mu = \mu (j)$.        
         Aus der Eigenschaft 1) haben wir $ a_{j_0 +2j+2}=1$ und
         $$
         q_{j_0 +2j+2} = q_{j_0 +2j+1} + q_{j_0 +2j}.
        $$
      Nun gilt
      $$
      s_{\mu+3}= q_{j_0 +2j+2} =q_{j_0 +2j+1}+ q_{j_0 +2j}
      < s_{\mu+2}+ s_{\mu+1}\le
      b_{\mu+3}s_{\mu+2}+ s_{\mu+1}=s_{\mu+3}.
      $$
      Es ist unm\"{o}glich.
      
        \noindent
   {\bf Teilfall 2.2.4.}
   Es gibt unendlich viele $j$ mit
   $$
   q_{j_0 +2j} = s_\mu, \,\,\,\,\,
      q_{j_0 +2j+2} = s_{\kappa},\,\,\,\,\,\,\,\, \kappa \ge\mu+4.
      $$
      Dann gibt es zwei sukzessive $s_{l-1}, s_l$  entweder in dem Intervall
      $(q_{j_0+2j}, q_{j_0+2j+1}]$,  oder  in dem Intervall $[q_{j_0+2j+1}, q_{j_0+2j+2})$.
      In beiden F\"{a}lle $ q_{i-1}<s_{l-1}< s_{l}<q_i$, wegen
Hilfssatz 3 ({\bf ii}). Hilfssatz 2    ({\bf ii}) gibt
entweder $\lambda_{l-1}(\beta) \ge K$ oder  $\lambda_{l}(\beta) \ge K$.
Schlie\ss lich  nach 
der Folgerung 2 
hat man  (\ref{ennq}) mit $m = l-1\,\,\text{ oder
}\,\, l$.

Also im allen F\"{a}llen entweder $\alpha\pm \beta \in \mathbb{Z}$, oder (\ref{enn}) oder (\ref{ennq}) gilt 
f\"{u}r beliebige gro\ss e Werten von Indizes. Satz 1 ist bewiesen.$\Box$

   {\bf 4. Beweis des Satzes 2.}
   
   Zu jedem $\varepsilon >0$ 
  gibt es ganzen Zahlen $U,V$ mit
   $$
   \left|V+\frac{U}{\tau} -
  \sqrt{\tau}\right|<\varepsilon.
   $$
   Definieren wir die Folge $X_n$ als
   $$
   X_0 = U,\,\,\,\,\, X_1 = V,\,\,\,\,\,\ X_{n+1}= X_n +X_{n-1}.
   $$
   Dann gilt
   $$
   X_n = 
   A\tau^n + B(-\tau)^{-n}
   \,\,\,\,\,
   \text{mit}
   \,\,\,\,\,
   A = \frac{\tau V+U}{\tau+2}.
   $$
   Wegen $ A>0$ gibt es $ k\in \mathbb{Z}_+$ mit
   $X_{k-1}, X_k\ge 1$.
   Sei
   $$
   \frac{X_{k-1}}{X_k} = [0; b_w,....,b_1],\,\,\,\,\, b_j \in \mathbb{Z}_+
   ,$$
    und definieren wir
    $$
    \varphi = [0,b_1,...,b_w,  \overline{1}].
    $$
    Dann
    \begin{equation}\label{pii}
    |\psi_\tau (t) - \psi_\varphi (t)|\le
    (K+\varepsilon)
    \min( \psi_\tau (t), \psi_\varphi (t) )
    \end{equation}
    f\"{u}r  alle hinreichende gro\ss e $t$.
    
    Tats\"{a}chlich,
    f\"{u}r $\tau$   ist $ q_n =F_n$ eine Fibonacci-Zahl und
    $$
    q_n = F_n= \frac{1}{\sqrt{5}}( \tau^n - (-\tau)^{-n}).
    $$
    Das Nenner $s_n$ des Naherungsbruchs $r_n/s_n$ 
    f\"{u}r $ \varphi$  ist
    $$
    s_n = X_{n-n_0},\,\,\,\,\,\, n_0={k-w} 
    $$
    f\"{u}r $n$ gro\ss \, genug. Bemerken wir, dass
    $$
    \frac{X_n}{F_n} 
    \sim A\sqrt{5} = U+V\tau^{-1},
    \,\,\,\,\,\,
       \left|A\sqrt{5} -
   \sqrt{\tau} \right|<\varepsilon .
    $$
    Dann
    $$q_n=F_{n} < X_n= s_{n+n_0} < F_{n+1}=q_{n+1}
    $$
    und
    $$
    \psi_\tau (t) =\xi_n = ||q_n\tau ||  \sim \frac{1}{\sqrt{5}F_n},\,\,\,\,\,\,
    F_n \le t < F_{n+1},
    $$
   $$
    \psi_\varphi (t) =\eta_n = ||s_{n+n_0}\varphi ||  \sim \frac{1}{\sqrt{5}X_n},\,\,\,\,\,\,
    X_n \le t < X_{n+1},
    $$
    $$
     \min( \psi_\tau (t), \psi_\varphi (t) )
     =
     \begin{cases}
     \xi_n,\,\,\, F_n \le t < X_n\cr
     \eta_n,\,\,\, x_n\le t < F_{n+1}
     \end{cases}.
    $$
    Sodass
    $$\frac{\xi_{n-1}}{\eta_n}
    \sim
    \frac{F_n}{X_{n-1}}  = \sqrt{\tau}+ O(\varepsilon)
   ,
    \,\,\,\,\,\,\,
    \frac{\eta_n}{\xi_n}    \sim
    \frac{X_n}{F_n}  = \sqrt{\tau}+ O(\varepsilon)
    $$
    und (\ref{pii}) ist bewiesen.$\Box$

   \end{document}